\title{ ~~\\ Convoluted convolved Fibonacci numbers}
\author{Pieter Moree}
\documentclass[12pt]{article}
\usepackage{amssymb, latexsym, amsfonts}
\def\@ptsize{2}
\newtheorem{Thm}{Theorem}
\newtheorem{Q}{Questions}
\newtheorem{Conj}{Conjecture}
\newtheorem{Lem}{Lemma}

\newtheorem{Cor}{Corollary}

\newtheorem{Prop}{Proposition}
\newcommand{\qed}{\hfill $\Box$}

\begin{document}
\date{}
\maketitle
{\def\thefootnote{}
}
\begin{abstract}
\noindent The convolved Fibonacci numbers $F_j^{(r)}$ are defined
by $(1-x-x^2)^{-r}=\sum_{j\ge 0}F_{j+1}^{(r)}x^j$. In this
note we consider some
related numbers that can be expressed in terms
of convolved Fibonacci numbers. These numbers appear in the numerical
evaluation of a constant arising in the study of the average density of
elements in a finite field having order congruent to $a({\rm mod~}d)$.
We derive a formula expressing these
numbers in terms of ordinary Fibonacci and Lucas numbers. The
non-negativity of these numbers can 
be inferred from
Witt's dimension formula for free Lie algebras.\\
\indent This note is a case study of the transform
${1\over n}\sum_{d|n}\mu(d)f(z^d)^{n/d}$ (with $f$ any formal series), which 
was introduced
and studied in the companion paper by Moree.
\end{abstract}
\section{Introduction}
Let $\{F_n\}_{n=0}^{\infty}=\{0,1,1,2,3,5,\dots\}$ be the sequence
of Fibonacci numbers and 
 $\{L_n\}_{n=0}^{\infty}=\{2,1,3,4,7,11,\dots\}$ the sequence of Lucas numbers.
It is well-known and easy to derive that
for $|z|<(\sqrt{5}-1)/2$, we have
$(1-z-z^2)^{-1}=\sum_{j=0}^{\infty}F_{j+1}z^j$. For any
real number $r$ the {\it convolved
Fibonacci numbers} are defined by
\begin{equation}
\label{fdefinitie}
{1\over (1-z-z^2)^r}=\sum_{j=0}^{\infty}F_{j+1}^{(r)}z^j.
\end{equation}
The Taylor series in (\ref{fdefinitie}) converges for all $z\in \mathbb C$
with $|z|< (\sqrt{5}-1)/2$. In the remainder of this note it is assumed that $r$ is
a positive integer.
Note that $F_{m+1}^{(r)}=\sum_{j_1+\cdots+j_r=m}F_{j_1+1}F_{j_2+1}\cdots F_{j_r+1}$,
where the sum is over all $j_1,\dots,j_r$ with $j_t\ge 0$ for $1\le t\le r$. We also
have $F_{m+1}^{(r)}=\sum_{j=0}^m F_{j+1}F_{m-j+1}^{(r-1)}$.\\
\indent The earliest reference to convolved 
Fibonacci numbers the author is aware of is a book by Riordan \cite{R}, who 
proposed as an exercise (at p. 89) to show that
$$F_{j+1}^{(r)}=\sum_{v=0}^{r}\left({r+j-v-1\atop j-v}\right)
\left({j-v\atop v}\right).$$
Convolved Fibonacci numbers have been studied in several papers,
for some references see, e.g., Sloane \cite{S}.
In 
Section \ref{reconsider} we give a formula expressing the convolved Fibonacci numbers
in terms of Fibonacci- and Lucas numbers. To the knowledge of the author it has
only been shown previously (by Hoggatt and 
Bicknell-Johnson \cite{HB} who used a different method) that this holds for $F_{j+1}^{(2)}$.\\
\indent In this note our main interest is in numbers $G_{j+1}^{(r)}$ and $H_{j+1}^{(r)}$ analogous
to the convolved Fibonacci numbers, which we name 
{\it convoluted convolved Fibonacci numbers}, respectively 
{\it sign twisted convoluted convolved Fibonacci numbers}. 
Given a formal series $f(z)\in \mathbb C[[z]]$, we define its {\it Witt transform} as
\begin{equation}
\label{definitie}
{\cal W}_f^{(r)}(z)={1\over r}\sum_{d|r}\mu(d)f(z^d)^{r\over d}=\sum_{j=0}^{\infty}m_f(j,r)z^j,
\end{equation}
where as usual $\mu$ denotes the M\"obius function.\\ 
\indent For every integer $r\ge 1$ we put 
$$G_{j+1}^{(r)}=m_f(j,r){\rm ~with~}f={1\over 1-z-z^2}.$$ 
Similarly we put
\begin{equation}
\label{definitie2}
H_{j+1}^{(r)}=(-1)^rm_f(j,r){\rm ~with~}f={-1\over 1-z-z^2}.
\end{equation}
On comparing (\ref{definitie2}) with (\ref{fdefinitie}) one sees that
\begin{equation}
\label{diri}
G_{j+1}^{(r)}={1\over r}\sum_{d|{\rm gcd}(r,j)}\mu(d)F_{{j\over d}+1}^{({r\over d})}{\rm ~and~}
H_{j+1}^{(r)}={(-1)^r\over r}\sum_{d|{\rm gcd}(r,j)}\mu(d)(-1)^{r\over d}
F_{{j\over d}+1}^{({r\over d})}.
\end{equation}
In Tables 1, 2  and 3 below some values of convolved,
convoluted convolved, respectively sign twisted convoluted convolved Fibonacci numbers are 
provided. 
The purpose of this paper is to investigate
the properties of these numbers. The next section gives a motivation
for studying these numbers.

\section{Evaluation of a constant}
Let $g$ be an integer and $p$ a prime not dividing $g$. Then by ord$_g(p)$ we denote
the smallest positive integer $k$ such that $g^k\equiv 1({\rm mod~}p)$. Let $d\ge 1$ be
an integer. It can be shown that the set of primes $p$ for which ord$_g(p)$ is divisible
by $d$ has a density and this density can be explicitly computed. It is easy to see
that the primes $p$ 
for which ord$_g(p)$ is even are the primes that divide some term of the sequence
$\{g^r+1\}_{r=0}^{\infty}$. A related, but 
much less studied, question is whether given integers $a$ and $d$
the set of primes $p$ for which ord$_g(p)\equiv a({\rm mod~}d)$ has a density. Presently
this more difficult problem can only 
be resolved under assumption of the Generalized Riemann Hypothesis, see 
, e.g., Moree \cite{M}. In the
explicit evaluation of this density and also that of its average value (where one
averages over $g$) \cite{PNew}, the following constant appears:
$$B_{\chi}=\prod_p\left(1+{[\chi(p)-1]p\over [p^2-\chi(p)](p-1)}\right),$$
where $\chi$ is a 
Dirichlet character and the product is over all primes $p$. Recall that the Dirichlet L-series for $\chi^k$, 
 $L(s,\chi^k)$, is defined, for Re$(s)>1$, by $\sum_{n=1}^{\infty}\chi^k(n)n^{-s}$. It satisfies
 the Euler product
 $$L(s,\chi^k)=\prod_p {1\over 1-{\chi^k(p)p^{-s}}}.$$
 We similarly define $L(s,-\chi^k)=\prod_p (1+\chi(p)p^{-s})^{-1}$. 
 The Artin constant, which appears in many problems involving the
 multiplicative order, is defined by
 $$A=\prod_p \left(1-{1\over p(p-1)}\right)=0.3739558136\dots$$
 The following result, which
follows from Theorem \ref{een} (with $f(z)=-z^3/(1-z-z^2)$) and some convergence arguments, 
expresses the constant $B_{\chi}$ in terms of Dirichlet L-series. Since Dirichlet
L-series in integer values are easily evaluated with very high decimal precision this
result allows one to evaluate $B_{\chi}$ with high decimal precision.
\begin{Thm} 
\label{Laatste}
{\rm 1)} We have
$$B_{\chi}=A{L(2,\chi)\over L(3,-\chi)}
\prod_{r=1}^{\infty}\prod_{j=3r+1}^{\infty}L(j,(-\chi)^r)^{-e(j,r)},$$
where $e(j,r)=G_{j-3r+1}^{(r)}$.\\
{\rm 2)} We have
$$B_{\chi}=A{L(2,\chi)L(3,\chi)\over L(6,\chi^2)}
\prod_{r=1}^{\infty}\prod_{j=3r+1}^{\infty}L(j,\chi^r)^{f(j,r)},$$
where $f(j,r)=(-1)^{r-1}H_{j-3r+1}^{(r)}$.
\end{Thm} 
{\it Proof}. Moree \cite{PNew} proved part 2, and a variation of his proof yields part 1. \qed\\ 

\noindent In the next
section it is deduced (Proposition \ref{propje})
that the numbers $e(j,r)$ appearing in the former double product are actually
positive integers and that the $f(j,r)$ are
non-zero integers that satisfy sgn$(f(j,r))=(-1)^{r-1}$. The proof makes
use of properties of the Witt transform that was introduced in \cite{PFibo}.

\section{Some properties of the Witt transform}
We recall some of the properties of the Witt transform 
(as defined by (\ref{definitie})) and deduce 
consequences for the (sign twisted) convoluted convolved Fibonacci numbers.
\begin{Thm} {\rm \cite{PNew}}.
\label{een}
Suppose that $f(z)\in \mathbb Z[[z]]$.
Then, as formal power series in $y$ and $z$, we have
$$1-yf(z)=\prod_{j=0}^{\infty}\prod_{r=1}^{\infty}(1-z^jy^r)^{m_f(j,r)}.$$
Moreover, the numbers $m_f(j,r)$ are integers. If
$$1-yf(z)=\prod_{j=0}^{\infty}\prod_{r=1}^{\infty}(1-z^jy^r)^{n(j,r)},$$
for some numbers $n(j,r)$, then $n(j,r)=m_f(j,r)$.
\end{Thm}
For certain choices of $f$ identities as above arise in the theory of Lie algebras, see, e.g., Kang 
and Kim \cite{KK}.
In this theory they go by the name of denominator identities.
\begin{Thm} 
\label{witttransform}
{\rm \cite{PFibo}}. Let $r\in \mathbb Z_{\ge 1}$ and $f(z)\in \mathbb Z [[z]]$. Write
$f(z)=\sum_j a_jz^j$.\\
{\rm 1)} We have
$$(-1)^r {\cal W}^{(r)}_{-f}(z)=\cases{{\cal W}^{(r)}_f(z)+{\cal
W}^{(r/2)}_f(z^2), &if $r\equiv 2({\rm mod~}4)$;\cr
{\cal W}^{(r)}_f(z), &otherwise.}$$
{\rm 2)} If $f(z)\in \mathbb Z [[z]]$, then so is ${\cal W}_f^{(r)}(z)$.\\
{\rm 3)} If $f(z)\in \mathbb Z_{\ge 0}[[z]]$, then so are ${\cal W}_f^{(r)}(z)$ and
$(-1)^r{\cal W}_{-f}^{(r)}(z)$.\\
\indent Suppose 
that $\{a_j\}_{j=0}^{\infty}$ is a non-decreasing sequence with $a_1\ge 1$.\\
{\rm 4)} Then $m_f(j,r)\ge 1$ and $(-1)^rm_{-f}(j,r)\ge 1$ for $j\ge 1$.\\
{\rm 5)} The
sequences $\{m_f(j,r)\}_{j=0}^{\infty}$ and $\{(-1)^rm_{-f}(j,r)\}_{j=0}^{\infty}$ 
are both non-decreasing.
\end{Thm}
In Moree \cite{PFibo} several further properties regarding monotonicity in both the $j$ and $r$ direction
are established that apply to both $G_j^{(r)}$ and $H_j^{(r)}$. It turns out that slightly stronger
results in this direction for these sequences can be established on using Theorem \ref{laatste} below.

\subsection{Consequences for $G_j^{(r)}$ and $H_j^{(r)}$}
Since clearly $F_{j+1}^{(r)} \in \mathbb Z$ we infer from (\ref{diri}) that
$rG_{j+1}^{(r)},rH_{j+1}^{(r)}\in \mathbb Z$. More is true, however:
\begin{Prop} 
\label{propje}
Let $j,r\ge 1$ be integers. Then\\
{\rm 1)} $G_{j}^{(r)}$ and $H_{j}^{(r)}$ are non-negative integers.\\ 
{\rm 2)} When $j\ge 2$, then $G_{j}^{(r)}\ge 1$ and $H_{j}^{(r)}\ge 1$.\\
{\rm 3)} We have
$$H_j^{(r)}=\cases{G_j^{(r)}+G_{j+1\over 2}^{(r/2)}, &if $r\equiv 2({\rm mod~}4)$ and $j$ is odd;\cr
G_j^{(r)}, &otherwise.}$$
{\rm 4)} The sequences $\{G_j^{(r)}\}_{j=1}^{\infty}$ and $\{H_j^{(r)}\}_{j=1}^{\infty}$ are
non-decreasing.
\end{Prop}
The proof easily follows from Theorem \ref{witttransform}.

\section{Convolved Fibonacci numbers reconsidered}
\label{reconsider}
We show that the convolved Fibonacci numbers can be expressed in terms of Fibonacci and
Lucas numbers. 
\begin{Thm} 
\label{fibolucas}

Let $j\ge 0$ and $r\ge 1$. We have 
$$F_{j+1}^{(r)}=\sum_{k=0\atop r+k\equiv 0({\rm mod~}2)}^{r-1}\left({r+k-1\atop k}\right)
\left({r-k+j-1\atop j}\right)
{L_{r-k+j}\over 5^{(k+r)/2}}+$$
$$\sum_{k=0\atop r+k\equiv 1({\rm mod~}2)}^{r-1}\left({r+k-1\atop k}\right)
\left({r-k+j-1\atop j}\right)
{F_{r-k+j}\over 5^{(k+r-1)/2}}.$$
In particular,
$5F_{j+1}^{(2)}=(j+1)L_{j+2}+2F_{j+1}$ and
$$50F_{j+1}^{(3)}=5(j+1)(j+2)F_{j+3}+6(j+1)L_{j+2}+12F_{j+1}.$$
\end{Thm}
{\it Proof}. Suppose that $\alpha,\beta\in \mathbb C$ with $\alpha \beta\ne 0$ and
$\alpha \ne \beta$. Then it is not difficult to show that we have the following
partial fraction decomposition:
$$(1-\alpha z)^{-r}(1-\beta z)^{-r}=$$
$$\sum_{k=0}^{r-1}\left({-r\atop k}\right)
{\alpha^r\beta^k\over (\alpha-\beta)^{r+k}}(1-\alpha z)^{k-r}+
\sum_{k=0}^{r-1}\left({-r\atop k}\right)
{\beta^r\alpha^k\over (\beta-\alpha)^{r+k}}(1-\beta z)^{k-r},$$
where $({-r\atop k})=1$ if $k=0$ and $({-r\atop k})=r(r-1)\cdots (r-k+1)/k!$ otherwise.
Using the Taylor expansion (with $t$ a real number)
$$(1-z)^t=\sum_{j=0}^{\infty}(-1)^j\left({t\atop j}\right)z^j,$$
we infer that $(1-\alpha z)^{-r}(1-\beta z)^{-r}=\sum_{j=0}^{\infty}\gamma(j)z^j$,
where 
$$\gamma(j)=\sum_{k=0}^{r-1}\left({-r\atop k}\right)
{\alpha^r\beta^k\over (\alpha-\beta)^{r+k}}(-1)^j\left({k-r\atop j}\right)\alpha^j+$$
$$\sum_{k=0}^{r-1}\left({-r\atop k}\right){\beta^r\alpha^k\over (\beta-\alpha)^{r+k}}(-1)^j\left({k-r\atop j}\right)\beta^j.$$
Note that $1-z-z^2=(1-\alpha z)(1-\beta z)$ with $\alpha=(1+\sqrt{5})/2$
and $\beta=(1-\sqrt{5})/2$. On substituting these values
of $\alpha$ and $\beta$ and using that
$\alpha-\beta=\sqrt{5}$, $\alpha \beta = -1$, $L_n=\alpha^n + \beta^n$ and
$F_n=(\alpha^n-\beta^n)/\sqrt{5}$, we find that
$$F_{j+1}^{(r)}=\sum_{k=0\atop r+k\equiv 0({\rm mod~}2)}^{r-1}(-1)^k\left({-r\atop k}\right)(-1)^j
\left({k-r\atop j}\right)
{L_{r-k+j}\over 5^{(k+r)/2}}+$$
$$\sum_{k=0\atop r+k\equiv 1({\rm mod~}2)}^{r-1}(-1)^k\left({-r\atop k}\right)(-1)^j
\left({k-r\atop j}\right)
{F_{r-k+j}\over 5^{(k+r-1)/2}}.$$
On noting that $(-1)^k({-r\atop k})=({r+k-1\atop k})$ and
$(-1)^j({k-r\atop j})=({r-k+j-1\atop j})$, the proof is completed. \qed\\

\noindent Let $r\ge 1$ be fixed. From the latter theorem one easily deduces the asymptotic
behaviour of $F_{j+1}^{(r)}$ considered as a function of $j$.
\begin{Prop}
\label{aso}
 Let $r\ge 2$ be fixed. Let $[x]$ denote the integer part of $x$. 
 Let $\alpha=(1+\sqrt{5})/2$. We have
 $F_{j+1}^{(r)}=g(r)j^{r-1}\alpha^j+O_r(j^{r-2}\alpha^j)$, as $j$ tends to infinity,
 where the implicit error term depends at most on $r$ and
 $g(r)=\alpha^r 5^{-[r/2]}/(r-1)!.$
 \end{Prop}
 
\section{The numbers $H_{j+1}^{(r)}$ for fixed $r$}

\noindent In this and the next section we consider the numbers  $H_{j+1}^{(r)}$ for fixed $r$, respectively
for fixed $j$. Very similar results can of course be obtained for the convoluted convolved Fibonacci
numbers $G_{j+1}^{(r)}$.\\
\indent For small fixed $r$ we can use Theorem \ref{fibolucas} in combination with (\ref{diri})
to explicitly express $H_{j+1}^{(r)}$ in terms of Fibonacci- and Lucas numbers. In doing so
it is convenient to work with the characteristic function $\chi$ of the integers, which is
defined by $\chi(r)=1$ if $r$ is an integer and $\chi(r)=0$ otherwise. We demonstrate the
procedure for $r=2$ and $r=3$. By (\ref{diri}) we find
$2H_{j+1}^{(2)}=F_{j+1}^{(2)}+F_{{j\over 2}+1}\chi({j\over 2})$ and
$3H_{j+1}^{(3)}=F_{j+1}^{(3)}-F_{{j\over 3}+1}\chi({j\over 3})$. By Theorem \ref{fibolucas} it
then follows for example that
$$150H_{j+1}^{(3)}=5(j+1)(j+2)F_{j+3}+6(j+1)L_{j+2}+12F_{j+1}-50F_{{j\over 3}+1}\chi({j\over 3}).$$
The asymptotic behaviour, for $r$ fixed and $j$ tending to infinity can be directly inferred
from (\ref{diri}) and Proposition \ref{aso}. 
\begin{Prop} With the same notation and assumptions as in Proposition {\rm \ref{aso}} we
have
 $H_{j+1}^{(r)}=g(r)j^{r-1}\alpha^j/r+O_r(j^{r-2}\alpha^j).$
\end{Prop}

\section{The numbers $H_{j+1}^{(r)}$ for fixed $j$}
In this section we investigate the numbers $H_{j+1}^{(r)}$ for fixed $j$. We first
investigate this question for the convolved Fibonacci numbers.\\
\indent The coefficient $F_{j+1}^{(r)}$ of $z^j$ in $(1-z-z^2)^{-r}$ is equal to the
coefficient of $z^j$ in $(1+F_2z+F_3z^3+\cdots+F_{j+1}z^j)^r$. By the
multinomial theorem we then find
\begin{equation}
\label{multinomial}
F_{j+1}^{(r)}=\sum_{\sum_{k=1}^j kn_k=j}\left({r\atop n_1,\cdots,n_j}\right)F_2^{n_1}
\cdots F_{j+1}^{n_j},
\end{equation}
where the multinomial coefficient is defined by
$$\left({r\atop m_1,\cdots,m_s}\right)=
{r!\over m_1!m_2!\cdots m_s!(r-m_1-\cdots-m_s)!}$$ and
$m_k\ge 0$ for $1\le k\le s$.\\

\noindent {\bf Example}. We have
\begin{eqnarray}
F_5^{(r)}&=&\left({r\atop 4}\right)+2\left({r\atop 2,1}\right)+4\left({r\atop 2}\right)
+3\left({r\atop 1,1}\right)+5\left({r\atop 1}\right)\nonumber\cr
&=&{7\over 4}r+{59\over 24}r^2+{3\over 4}r^3+{r^4\over 24}.\nonumber
\end{eqnarray}
This gives an explicit description of the sequence $\{F_5^{(r)}\}_{r=1}^{\infty}$ which is
sequence A006504 of \cite{S}.\\

\noindent The sequence $\{({r\atop m_1,\cdots,m_k})\}_{r=0}^{\infty}$ is a polynomial
sequence where the degree of the polynomial is $m_1+\cdots+m_k$. It follows from this
and (\ref{multinomial}) that
$\{F_{j+1}^{(r)}\}_{r=0}^{\infty}$ is a polynomial sequence of degree
$\max\{n_1+\cdots+n_j|\sum_{k=1}^j kn_j=j\}=j$. The leading term of this polynomial in $r$
is due to the multinomial term having $n_1=j$ 
and $n_t=0$ for $2\le t\le j$. All other terms in (\ref{multinomial}) are of lower degree. We thus
infer that $F_{j+1}^{(r)}=r^j/j!+O_j(r^{j-1}),~r\rightarrow \infty$. 
We leave it to the reader to make this more precise by showing that the coefficient of
$r^{j-1}$ is $3/(2(j-2)!)$.
If $n_1,\dots,n_j$ 
satisfy $\max\{n_1+\cdots+n_j|\sum_{k=1}^j kn_j=j\}=j$, then 
$j!/(n_1!\dots n_j!)$ is an integral multiple of a multinomial coefficient and hence
an integer. We thus infer that $j!F_{j+1}^{(r)}$ is a monic polynomial
in $\mathbb Z[r]$ of degree $j$. Note that the constant term of this polynomial is zero.
To sum up, we have obtained:
\begin{Thm}
\label{voorlaatste}
Let $j,r\ge 1$ be integers. There is a polynomial
$$A(j,r)=r^j+{3\over 2}j(j-1)r^{j-1}+\cdots \in \mathbb Z[r]$$ 
with $A(j,0)=0$ such that
$F_{j+1}^{(r)}=A(j,r)/j!$.
\end{Thm}
Using this result, the following regarding the sign twisted convoluted convolved Fibonacci numbers
can be established.
\begin{Thm}
\label{twee}
Let $\chi(r)=1$ if $r$ is an integer and $\chi(r)=0$ otherwise. We have
$$H_1^{(r)}=\cases{1, &if $r\le 2$;\cr 0, &otherwise,}$$
furthermore $H_2^{(r)}=1$. We have
$$2H_3^{(r)}={3}+{r}-{(-1)^{r/2}}\chi({r\over 2}){\rm ~and~}
6H_4^{(r)}={8}+{9}r+{r^2}-{2}\chi({r\over 3}).$$
Also we have
$$24H_5^{(r)}={42}+{59}r+{18}r^2+{r^3}-(18+{3r})
(-1)^{r\over 2}\chi({r\over 2}){\rm ~and~}$$
$$120H_6^{(r)}=264+450r+215r^2+{30r^3}+{r^4}-
{24}\chi({r\over 5}).$$
In general we have
$$H_{j+1}^{(r)}=\sum_{d|j,~2\nmid d}\mu(d)\chi({r\over d}){A({j\over d},{r\over d})\over 
r(j/d)!}+
\sum_{d|j,~2|d}\mu(d)(-1)^{r/2}\chi({r\over d}){A({j\over d},{r\over d})\over 
r(j/d)!}.$$
Let $j\ge 3$ be fixed. As $r$ tends to infinity we have
$$H_{j+1}^{(r)}={r^{j-1}\over j!}+{3r^{j-2}\over 2(j-2)!}+O_j(r^{j-2}).$$
\end{Thm}
{\it Proof}. Using that $\sum_{d|n}\mu(d)=0$ if $n>1$, it  is easy to check that
$$H_1^{(r)}={(-1)^r\over r}\sum_{d|r}\mu(r)(-1)^{r/d-1}
=\cases{1, &if $r\le 2$;\cr
0, &if $r>2$.}$$
The remaining assertions can be all derived from (\ref{diri}), (\ref{multinomial}) 
and Theorem \ref{voorlaatste}. \qed

\section{Monotonicity}
Inspection of the tables below suggests monotonicity properties of
$F_j^{(r)},G_j^{(r)}$ and $H_j^{(r)}$ to hold true. 
\begin{Prop} $~$\\
{\rm 1)} Let $j\ge 2$. Then $\{F_j^{(r)}\}_{r=1}^{\infty}$
is a strictly increasing sequence.\\
{\rm 2)} Let $r\ge 2$. Then $\{F_j^{(r)}\}_{j=1}^{\infty}$
is a strictly increasing sequence.
\end{Prop}
The proof of this is easy. For the proof of part 2 one can make use of the following simple observation.
\begin{Lem}
\label{ergflauw}
Let $f(z)=\sum_j a(j)z^j\in \mathbb R[[z]]$ be a formal series. Then $f(z)$ is said to have $k$-nondecreasing
coefficients if $a(k)>0$ and $a(k)\le a(k+1)\le a(k+2)<\dots$. If $a(k)>0$
and $a(k)<a(k+1)<a(k+2)<\dots$, then $f$ is said to have $k$-increasing coefficients.\\
If $f,g$ are $k$-increasing, respectively $l$-nondecreasing, then $fg$ is $(k+l)$-increasing.\\
If $f$ is $k$-increasing and $g$ is $l$-nondecreasing, then $f+g$ is $\max(k,l)$-increasing.\\
If $f$ is $k$-increasing, then $\sum_{j\ge 1}b(j)f^j$ with $b(j)\ge 0$ and $b(1)>0$ is $k$-increasing.
\end{Lem} 
We conclude this paper by establishing the following result:
\begin{Thm} \label{monotonie} $~$\\
{\rm 1)} Let $j\ge 4$. Then $\{G_j^{(r)}\}_{r=1}^{\infty}$
is a strictly increasing sequence.\\
{\rm 2)} Let $r\ge 1$. Then $\{G_j^{(r)}\}_{j=2}^{\infty}$
is a strictly increasing sequence.\\
{\rm 3)} Let $j\ge 4$. Then $\{H_j^{(r)}\}_{r=1}^{\infty}$
is a strictly increasing sequence.\\
{\rm 4)} Let $r\ge 1$. Then $\{H_j^{(r)}\}_{j=2}^{\infty}$
is a strictly increasing sequence.
\end{Thm}
The proof rests on expressing the entries of the above sequences in terms of certain quantities occurring in
the theory of free Lie algebras and circular words 
(Theorem \ref{laatste}) and then invoke results on the monotonicity of
these quantities to establish the result.

\subsection{Circular words and Witt's dimension formula}
We will make use of an easy result on cyclic words. A word $a_1\cdots a_n$ is
called {\it circular} 
or {\it cyclic} if $a_1$ is regarded as following $a_n$, where
$a_1a_2\cdots a_n$,
$a_2\cdots a_na_1$
and all other cyclic shifts (rotations)
of $a_1a_2\cdots a_n$ are regarded as the same word. A circular
word of length $n$ may conceivably be given by repeating a segment of
$d$ letters $n/d$
times, with $d$ a divisor of $n$. Then one says the word is of {\it period} $d$.
Each word belongs
to an unique smallest period: the {\it minimal period}.\\
\indent Consider
circular words of length $n$ on an alphabet $x_1,\dots,x_r$ consisting of
$r$ letters.
The total number of ordinary words such that $x_i$ occurs $n_i$ times equals
$${n!\over n_1!\cdots n_r!},$$
where $n_1+\cdots+n_r=n$.
Let $M(n_1,\dots,n_r)$ denote
the number of circular words of length $n_1+\cdots+n_r=n$ and
minimal period $n$ such that the letter $x_i$ appears
exactly $n_i$ times.
This leads to the formula
\begin{equation}
\label{witty}
{n!\over n_1!\cdots n_r!}=\sum_{d|{\rm
gcd}(n_1,\dots,n_r)}{n\over d}M({n_1\over d},{n_2\over d},
\dots,{n_r\over d}).
\end{equation}
whence it follows by M\"obius inversion that
\begin{equation}
\label{basaal}
M(n_1,\dots,n_r)={1\over n}
\sum_{d|{\rm gcd}(n_1,\dots,n_r)}\mu(d){{n\over d}!\over {n_1\over
d}!\cdots {n_r\over d}!}.
\end{equation}
Note that $M(n_1,\dots,n_r)$ is totally symmetric in the variables $n_1,\dots,n_r$.
\indent The numbers $M(n_1,\dots,n_r)$ 
also occur in a classical result
in Lie theory, namely Witt's formula for the homogeneous subspaces of a
finitely generated
free Lie algebra $L$: if $H$ is the subspace of $L$ generated by all
homogeneous elements of
multidegree $(n_1,\dots,n_r)$, then dim$(H)=M(n_1,\dots,n_r)$, where
$n=n_1+\cdots+n_r$.\\ 
\indent In Theorem \ref{laatste} a variation of $M(n_1,\dots,n_r)$ appears.
\begin{Lem} {\rm \cite{PFibo}}.
\label{minnetje}
Let $r$ be a positive integer and let $n_1,\dots,n_r$ be non-negative
integers and put $n=n_1+\cdots+n_r$. Let
$$V_1(n_1,\dots,n_r)={(-1)^{n_1}\over
n}\sum_{d|{\rm gcd}(n_1,\dots,n_r)}\mu(d)(-1)^{n_1\over d}
{{n\over d}!\over {n_1\over d}!\cdots {n_r\over d}!}.$$
Then $$V_1(n_1,\dots,n_r)=$$
$$\cases{
M(n_1,\dots,n_r)+M({n_1\over 2},\dots,{n_r\over 2}), &if $n_1\equiv 2({\rm
mod~}4)$ and
$2|{\rm gcd}(n_1,\dots,n_r)$,;\cr
M(n_1,\dots,n_r) & otherwise.}$$
\end{Lem} 
The numbers $V_1(n_1,\dots,n_r)$ can also be interpreted as dimensions (in the context
of free Lie superalgebras), see, e.g., Petrogradsky \cite{P}.\\
\indent The numbers $M$ and $V_1$ enjoy certain monotonicity properties. 
\begin{Lem} 
\label{atlonglast} {\rm \cite{PFibo}}.
Let $r\ge 1$ and $n_1,\dots,n_r$ be non-negative numbers.\\
{\rm 1)}  The
sequence $\{M(m,n_1,\dots,n_r)\}_{m=0}^{\infty}$ is 
non-decreasing  if $n_1+\cdots+n_r\ge 1$ and
strictly
increasing if $n_1+\cdots+n_r\ge 3$.\\
{\rm 2)} The sequence 
$\{V_1(m,n_1,\dots,n_r\}_{m=0}^{\infty}$ is non-decreasing  if $n_1+\cdots+n_r\ge 1$
and strictly
increasing if $n_1+\cdots+n_r\ge 3$. 
\end{Lem}
\indent Using (\ref{witty}) one infers (on taking the logarithm of either side and
expanding it as a formal series) that
\begin{equation}
\label{fundamental}
{1-z_1-\cdots
-z_r}=\prod_{n_1,\dots,n_r=0}^{\infty}(1-z_1^{n_1}\cdots
z_r^{n_r})^{M(n_1,\cdots,n_r)},
\end{equation}
where $(n_1,\dots,n_r)=(0,\dots,0)$ is excluded in the product. From the latter identity it
follows that
$$
1+z_1-z_2-\cdots
-z_r=$$
$$\prod_{n_1,\dots,n_r=0\atop 2|n_1}^{\infty}(1-z_1^{n_1}\cdots
z_r^{n_r})^{M(n_1,\dots,n_r)}\prod_{n_1,\dots,n_r=0\atop 2\nmid n_1}^{\infty}\left(1-z_1^{2n_1}\cdots
z_r^{2n_r}\over 1-z_1^{n_1}\cdots z_r^{n_r}\right)^{M(n_1,\dots,n_r)},
$$
whence, by Lemma \ref{minnetje}, 
\begin{equation}
\label{fundamental2}
1+z_1-z_2-\cdots-z_r
=\prod_{n_1,\dots,n_r=0}^{\infty}(1-z_1^{n_1}\cdots
z_r^{n_r})^{(-1)^{n_1}V_1(n_1,\dots,n_r)}.
\end{equation}

\begin{Thm}
\label{laatste}
Let $r\ge 1$ and $j\ge 0$. We have
$$G_{j+1}^{(r)}=\sum_{k=0}^{[j/2]}M(r,k,j-2k) {\rm ~and~}H_{j+1}^{(r)}=\sum_{k=0}^{[j/2]}V_1(r,k,j-2k).$$
\end{Thm}
{\it Proof}. By Theorem \ref{een} and the definition of $G_j^{(r)}$ we infer that
$$1-{y\over 1-z-z^2}=\prod_{j=0}^{\infty}\prod_{r=1}^{\infty}(1-z^jy^r)^{G_{j+1}^{(r)}}.$$
The left hand side of the latter equality equals $(1-z-z^2-y)/(1-z-z^2)$. On invoking 
(\ref{fundamental}) with 
$z_1=z$, $z_2=z^2$ and $z_3=y$ the claim
regarding $G_{j}^{(r)}$ follows from the uniqueness assertion in Theorem \ref{een}.\\ 
\indent The proof of the identity for $H_{j+1}^{(r)}$ is similar, but makes use of
identity (\ref{fundamental2}) instead of (\ref{fundamental}). \qed\\

\noindent {\it Proof of Theorem} \ref{monotonie}. 1) For $j\ge 5$, $k+j-2k\ge j-[j/2]\ge 3$ and hence
each of the terms $M(r,k,j-2k)$ with $0\le k\le [j/2]$ is strictly increasing in $r$ by Lemma \ref{atlonglast}.
For $3\le j\le 4$, by Lemma \ref{atlonglast} again, all terms $M(r,k,j-2k)$ with $0\le k\le [j/2]$ are
non-decreasing in $r$ and at least one of them is strictly increasing. The result now follow
by Theorem \ref{laatste}.\\
2) In the proof of part 1 replace the letter `$M$' by `$V_1$'.\\
3) For $r=1$ we have $G_{j+!}^{(1)}=F_{j+1}$ and the
the result is obvious. For $r=2$ 
each of the terms $M(r,k,j-2k)$ with $0\le k\le [j/2]$ is non-decreasing in $j$. For $j\ge 2$
one of these is strictly increasing. Since in 
addition $G_2^{(2)}<G_3^{(2)}$ the result follows for $r=2$.
For $r\ge 3$ each of the terms $M(r,k,j-2k)$ with $0\le k\le [j/2]$ is strictly increasing in $j$.
The result now follows by Theorem \ref{laatste}.\\
4) In the proof of part 3 replace the letter `$G$' by `$H$ and `$M$' by `$V_1$'. \qed

\vfil\eject
\hfil\break
\section{Tables}

\centerline{{\bf Table 1:} Convolved Fibonacci numbers $F_j^{(r)}$}
\medskip
\begin{center}
\begin{tabular}{|c|c|c|c|c|c|c|c|c|c|c|c|}\hline
$r\backslash j$&1&2&3&4&5&6&7&8&9&10&11\\ \hline
\hline\hline
1&1&1&2&3&5&8&13&21&34&55&89\\ \hline
2&1&2&5&10&20&38&71&130&235&420&744\\ \hline
3&1&3&9&22&51&111&233&474&942&1836&3522\\ \hline
4&1&4&14&40&105&256&594&1324&2860&6020&12402\\ \hline
5&1&5&20&65&190&511&1295&3130&7285&16435&36122\\ \hline
\end{tabular}  
\end{center}

\centerline{{\bf Table 2:} Convoluted convolved Fibonacci numbers $G_j^{(r)}$}
\medskip
\begin{center}
\begin{tabular}{|c|c|c|c|c|c|c|c|c|c|c|c|}\hline
$r\backslash j$&1&2&3&4&5&6&7&8&9&10&11\\ \hline
\hline\hline
1&1&1&2&3&5&8&13&21&34&55&89\\ \hline
2&0&1&2&5&9&19&34&65&115&210&368\\ \hline
3&0&1&3&7&17&37&77&158&314&611&1174\\ \hline
4&0&1&3&10&25&64&146&331&710&1505&3091\\ \hline
5&0&1&4&13&38&102&259&626&1457&3287&7224\\ \hline
6&0&1&4&16&51&154&418&1098&2726&6570&15308\\ \hline
7&0&1&5&20&70&222&654&1817&4815&12265&30217\\ \hline
8&0&1&5&24&89&309&967&2871&8043&21659&56123\\ \hline
9&0&1&6&28&115&418&1396&4367&12925&36542&99385\\ \hline
10&0&1&6&33&141&552&1946&6435&20001&59345&168760\\ \hline
\end{tabular}  
\end{center}

\centerline{{\bf Table 3:} Sign twisted convoluted convolved Fibonacci numbers $H_j^{(r)}$}
\medskip
\begin{center}
\begin{tabular}{|c|c|c|c|c|c|c|c|c|c|c|c|}\hline
$r\backslash j$&1&2&3&4&5&6&7&8&9&10&11\\ \hline
\hline\hline
1&1&1&2&3&5&8&13&21&34&55&89\\ \hline
2&1&1&3&5&11&19&37&65&120&210&376\\ \hline
3&0&1&3&7&17&37&77&158&314&611&1174\\ \hline
4&0&1&3&10&25&64&146&331&710&1505&3091\\ \hline
5&0&1&4&13&38&102&259&626&1457&3287&7224\\ \hline
6&0&1&5&16&54&154&425&1098&2743&6570&15345\\ \hline
7&0&1&5&20&70&222&654&1817&4815&12265&30217\\ \hline
8&0&1&5&24&89&309&967&2871&8043&21659&56123\\ \hline
9&0&1&6&28&115&418&1396&4367&12925&36542&99385\\ \hline
10&0&1&7&33&145&552&1959&6435&20039&59345&168862\\ \hline
\end{tabular}  
\end{center}

\hfil\break

\medskip\noindent {\footnotesize Korteweg-de Vries Institute,
Plantage Muidergracht 24, 1018 TV Amsterdam, The Netherlands.\\
e-mail: moree@science.uva.nl }

\end{document}